\documentclass{article}

\usepackage{amsmath}
\usepackage{amsthm}
\usepackage{amssymb}
\usepackage{latexsym}
\usepackage{enumerate}

\def\M{{\mathcal M}}

\def\Z{{\mathbb Z}}
\def\D{{\mathbb D}}
\def\A{{\mathbb A}}
\def\S{{\mathbb S}}

\begin{document}

\markboth{M.Mecchia and B.P.Zimmermann}{On minimal finite factor groups}

\title{On minimal finite factor groups of outer automorphism
groups of free groups}

\author{Mattia Mecchia and Bruno P. Zimmermann}

\date{}

\maketitle

\begin{abstract}
We prove that, for $n=3$ and 4, the minimal nonabelian finite
factor group of the outer automorphism group ${\rm Out} \, F_n$ of a free group
of rank $n$ is the linear group ${\rm PSL}_n(\Z_2)$ (conjecturally, this may
remain true for arbitrary rank $n > 2$). We also discuss some computational
results on low index subgroups of  ${\rm Aut} \, F_n$ and ${\rm Out} \, F_n$,
for $n = 3$ and 4, using presentations of these groups.
\end{abstract}

\maketitle

\bigskip \bigskip

\section{Introduction}

\medskip

It is shown in \cite{Z} that the minimal nontrivial finite quotient (nontrivial factor
group of smallest possible order) of the mapping class group  $\M_g$ of a closed
orientable surface of genus $g$ is the symplectic group ${\rm PSp}_{2g}(\Z_2)$,
for $g = 3$ and 4; in fact, this may remain true for arbitrary genus $g \ge 3$. Since
$\M_g$ is perfect for $g \ge 3$, such a minimal nontrivial finite quotient is a
nonabelian simple group.  There are canonical projections onto symplectic groups
$$\M_g  \to  {\rm Sp}_{2g}(\Bbb Z)  \to {\rm Sp}_{2g}(\Z_p)
\to  {\rm  PSp}_{2g}(\Z_p),$$ and the projective symplectic groups ${\rm
PSp}_{2g}(\Z_p)$ are simple if $p$ is prime and $g \ge 3$. It is a consequence of the
congruence subgroup property for the symplectic groups ${\rm Sp}_{2g}(\Bbb Z)$ that
their finite simple quotients are exactly the finite projective symplectic groups
${\rm PSp}_{2g}(\Z_p)$ (see \cite{Z}). But also for the mapping class groups $\M_g$,
all known finite quotients seem to be strongly connected to the symplectic
groups, and it would be interesting to know what other finite simple groups can
occur (see \cite{T} for a computational approach for genus two and three, and \cite{MR}
for comments on the congruence subgroup property).

\medskip

In the present note, we consider the outer automorphism group ${\rm Out} \,
F_n$ of a free group $F_n$ of rank $n$. It is well-known that ${\rm Out} \, F_2
\cong {\rm GL}_2(\Z) \cong \D_{12} *_{\D_4}\D_8$ (a free product with
amalgamation of two dihedral groups of orders 12 and 8), and we shall assume in
the following that $n \ge 3$; then  the abelianization of ${\rm Out} \, F_n$
has order two.  There is a canonical projection of ${\rm Out} \, F_n$ onto ${\rm
GL}_n(\Z)$, and we consider also the preimage of ${\rm SL}_n(\Z)$ in  ${\rm
Out} \, F_n$ which we denote by  ${\rm SOut} \, F_n$ (the unique subgroup of
index two of ${\rm Out} \, F_n$). It is well-known that ${\rm SOut} \, F_n$ is
a perfect group (see \cite{Ge} for a presentation), so the minimal nontrivial
quotient will be again a nonabelian simple group. There are projections
$${\rm SOut} \, F_n \to {\rm SL}_n(\Bbb Z)  \to {\rm SL}_n(\Z_p) \to  {\rm
PSL}_n(\Z_p),$$

and the finite linear groups ${\rm PSL}_n(\Z_p)$ are simple if $p$ is prime.  It is
a consequence of the congruence subgroup property for the linear group
${\rm SL}_n(\Z)$ that the finite simple quotients of ${\rm SL}_n(\Z)$ are exactly
the finite projective linear groups  ${\rm PSL}_n(\Z_p)$, $p$ prime (\cite{Z}).

\medskip

Our main result is the following:

\bigskip

{\bf Proposition.} {\sl For  $n=3$ and 4, the minimal nontrivial finite quotient
of ${\rm SOut} \, F_n$, and also the minimal nonabelian finite quotient of  ${\rm Out}
\, F_n$, is the linear group ${\rm PSL}_n(\Z_2)$.}

\bigskip

We note that ${\rm PSL}_3(\Z_2)  \cong {\rm PSL}_2(\Z_7)$ is the unique simple
group of order 168, and that ${\rm PSL}_4(\Z_2)$ is isomorphic to the alternating
group $\A_8$, of order 20160. Conjecturally, the Proposition  remains true for
arbitrary $n \ge 3$. Concerning other finite simple groups, it is shown in \cite{Gi} that
infinitely many alternating groups occur as quotients of ${\rm Out} \, F_n$.

\medskip

As a consequence of the Proposition we have also the following:

\bigskip

{\bf Corollary.} {\sl  The minimal index of a proper subgroup of  ${\rm SOut}
\, F_4$, and also of a proper subgroup of ${\rm Out} \, F_4$ different from
${\rm SOut} \, F_4$, is eight (the minimal index of a proper subgroup of
${\rm PSL}_4(\Z_2)  \cong \A_8$).}

\bigskip

For $n = 3$ this minimal index should be seven (the minimal index of a proper
subgroup of ${\rm PSL}_3(\Z_2)  \cong {\rm PSL}_2(\Z_7)$) but for the moment we
cannot exclude index six by the present methods; we verified this, however, by
computational methods (GAP), see section 3 for some comments.

\bigskip

\section{Proof of the Proposition and the Corollary}

\bigskip

We denote by  ${\rm Aut} \, F_n$ the automorphism group of the free group $F_n$
and by ${\rm SAut} \, F_n$ its subgroup of index two which is the preimage of
${\rm SL}_n(\Z)$ under the canonical projection of ${\rm Aut} \, F_n$  onto
${\rm GL}_n(\Z)$.  Fixing a free generating set of $F_n$, inversions and
permutations of generators generate a subgroup (Weyl group)
$W_n \cong (\Z_2)^n \rtimes \S_n$ of  ${\rm Aut} \, F_n$; let $SW_n$ denote $W_n \cap
\, {\rm SAut} \, F_n$, with $SW_n \cong (\Z_2)^{n-1} \rtimes \S_n$.

\medskip

We note that, by results in \cite{WZ}, $W_n$ is  the finite subgroup of maximal possible
order of both ${\rm Aut} \, F_n$ and ${\rm Out} \, F_n$, for $n \ge 3$, unique up to
conjugation if $n > 3$ (for $n=3$ there is one other subgroup of maximal possible
order 48).

\medskip

Let  $\Delta$ denote the central element of
$W_n$ inverting all generators; note that
$\Delta$ is in  $SW_n$ if and only if $n$ is even.

\medskip

The proof of the Proposition is based on the following:

\bigskip

{\bf Lemma} (\cite[Prop. 3.1]{BV}).  {\sl  Let $n \ge 3$ and $\phi$ be a
homomorphism from ${\rm SAut} \, F_n$ to a group $G$. If the restriction of
$\phi$ to $SW_n$ has nontrivial kernel $K$ then one of the following holds:

\begin{enumerate}[i)]

\item $n$ is even, $K = \langle \Delta \rangle$ and $\phi$ factors through
${\rm PSL}_n(\Z)$;

\item $K$ is the intersection of $SW_n$ with the subgroup $(\Z_2)^n$ of $W_n$
generated by all inversions and the image of $\phi$ is isomorphic
to ${\rm PSL}_n(\Z_2)$, or

\item $\phi$  is the trivial map.

\end{enumerate} }

\bigskip

{\it Proof of the Proposition.}  We consider the case $n = 3$ first.  Since
${\rm SOut} \, F_3$ is perfect, a  minimal nontrivial finite quotient of ${\rm
SOut} \, F_3$ is a nonabelian simple group. The only nonabelian simple group
with an order smaller than the order 168 of the linear group  ${\rm PSL}_3(
\Z_2)  \cong {\rm PSL}_2(\Z_7)$ is the alternating group $\A_5$, of order 60
(see \cite{C} for information about the finite simple groups).  Since the order 24
of $SW_3$ does not divide the order of $\A_5$, by the Lemma every
homomorphism from  ${\rm SOut} \, F_3$ to $\A_5$ is trivial, hence the minimal
possibility for a nontrivial finite quotient of ${\rm SOut} \, F_3$ is the
linear group ${\rm PSL}_3(\Z_2)  \cong {\rm PSL}_2(\Z_7)$, the unique simple
group of order 168.

\medskip

As for ${\rm Out} \, F_3$, if the finite nonabelian group $G$ is  a quotient of
${\rm Out} \, F_3$ then the image of
${\rm SOut} \, F_3$ has index one or two in $G$ and order at least 168.
Since ${\rm Out} \, F_3$ surjects onto ${\rm PSL}_3(\Z_2) =  {\rm PGL}_3( \Z_2)$,
of order 168, this is again the minimal possibility for $G$.

\bigskip

We come now to the proof of the Proposition for $n=4$. Suppose that $\phi: {\rm SOut}
\, F_4 \to G$ is a nontrivial homomorphism onto a finite simple group $G$ of order
less than the order 20160 of ${\rm PSL}_4( \Z_2) \cong \A_8$. The simple groups
of order less than 20160 are the following (see \cite{C}):

\begin{itemize}

\item the alternating groups $\A_d$ of degrees $d=5, 6$ or 7;

\item the linear groups ${\rm PSL}_2(\Z_p) = {\rm L}_2(p)$, for the primes $p=7$, 11, 13,
17, 19, 23, 29 and 31;

\item the linear groups  ${\rm PSL}_2(q) = {\rm L}_2(q)$, for the prime powers $q=8$, 9,
16, 25 and 27 (over the finite fields of the corresponding orders);

\item the linear group ${\rm PSL}_3( \Z_3) = {\rm L}_3(3)$;

\item the unitary group ${\rm PSU}_3( \Z_3) = {\rm U}_3(3)$;

\item the Mathieu group ${\rm M}_{11}$.

\end{itemize}

\medskip

Since the order of none of these groups is divided by the order $2^6 \cdot 3$
of $SW_4$, the restriction of $\phi$ to $SW_4$ has nontrivial kernel and
the Lemma applies. Since the cases ii) and iii) of the Lemma are
excluded by the hypotheses, we are necessarily in case i). Now also case i) may
be excluded by appealing to the fact that the finite quotients of ${\rm
PSL}_4(\Z)$ are exactly the groups  ${\rm PSL}_4(\Z_p)$, $p$ prime (\cite[Theorem
1]{Z}); however also the following more direct argument, in the spirit of the
previous ones, applies. In case i) of the Lemma, the kernel of the
restriction of $\phi$ to $SW_4$ is the subgroup of order two generated by the
central involution $\Delta$ of $SW_4$. Hence the factor group $SW_4/\langle
\Delta \rangle$, of order $2^5 \cdot 3$, embeds into $G$. However none of the
above groups has such a subgroup: considering orders again, the only candidates
which remain are ${\rm U}_3(3)$ and ${\rm L}_2(31)$, but both do not have a
subgroup isomorphic to $SW_4/\langle \Delta \rangle$.

\medskip

There is one other simple group of order 20160 not isomorphic to ${\rm
PSL}_4(\Z_2) \cong \A_8$,  the linear group ${\rm PSL}_3(4) = {\rm L}_3(4)$.
Considering the maximal subgroups of ${\rm L}_3(4)$, it is easy to see
that $G = {\rm L}_3(4)$ has no subgroup isomorphic to
$SW_4$, so the Lemma applies again leaving us with case i); finally, also
case i) is excluded since ${\rm L}_3(4)$ has no subgroup isomorphic to
$SW_4/\langle \Delta \rangle$ (or by appealing again to \cite[Theorem 1]{Z}).

\medskip

This completes the proof of the Proposition for ${\rm SOut} \, F_4$, and hence also
for ${\rm Out} \, F_4$.

\bigskip

{\it Proof of the Corollary.}  Both ${\rm SOut} \, F_4$ and ${\rm Out} \, F_4$
have subgroups of index eight since both admit surjections onto the alternating
group $\A_8  \cong {\rm PSL}_4(\Z_2)$.  By the proof of the Proposition, ${\rm
SOut} \, F_4$ does not admit a nontrivial homomorphism to an alternating groups
$\A_d$ of degree $d < 8$, hence ${\rm SOut} \, F_4$ has no proper subgroup of index
less than eight. Similarly, the same holds for ${\rm Out} \, F_4$ which does
not admit a nontrivial homomorphism to a symmetric group $\S_d$ of degree $d < 8$.

\bigskip

\section{Comments on computations}

\bigskip

We employed the low index subgroup procedure of the computer algebra system GAP
in order to find the smallest indices of proper subgroups of ${\rm SAut} \,
F_3$ and ${\rm SOut} \, F_3$. We used the 4-generator presentation of ${\rm
Aut} \, F_3$ in \cite[section 3.5, Corollary N1]{MaKS} (note that some of the
relations apply only for $n > 3$; see also \cite[section 7.3]{CM}, and \cite{MC}, \cite[section 6]{NN},
section 6] for presentations of ${\rm Out} \, F_3$), and created by GAP the
unique subgroups ${\rm SAut} \, F_3$ and  ${\rm SOut} \, F_3$ of index two. We
found that the three smallest indices of a proper subgroup of ${\rm SAut} \,
F_3$ and ${\rm SOut} \, F_3$  are  7, 8 and 13; the factor groups of the cores
of the corresponding subgroups (the largest normal subgroup contained in a
subgroup) are  ${\rm PSL}_3(\Z_2) = {\rm GL}_3(\Z_2)$ for indices 7 and 8, a
semidirect product  $(\Z_2)^3 \rtimes {\rm GL}_3(\Z_2)$ for index 8, and ${\rm
PSL}_3(\Z_3)$ for index 13.

\bigskip

Remark. Concerning the smallest indices of proper subgroups of ${\rm SL}_3(\Bbb
Z)$, index 8 is now missing and there remain only the indices 7 and 13 (see
\cite[section 3.5]{MaKS} for a presentation of ${\rm GL}_r(\Bbb Z)$). Computing the
abelianization of the core of the index 8 subgroup of  ${\rm SOut} \, F_3$ (a
normal subgroup of index 1344, with quotient $(\Z_2)^3 \rtimes {\rm
GL}_3(\Z_2)$), we found the free abelian group $\Bbb Z^{14}$;  so this gives an
explicit example of a finite index subgroup of ${\rm
SOut} \, F_3$ with infinite abelianization.  We note that the existence of such
subgroups was known since by \cite{MC}, ${\rm Out} \, F_3$ is virtually residually
torsion-free nilpotent; on the other hand, no such examples seem to be known
for rank $n > 3$ (see \cite{T} for the corresponding situation for mapping class
groups).

\bigskip

We then went on to compute that the three minimal indices
of subgroups of  ${\rm SAut} \, F_4$ are 8,  15 and 16; the factor groups
of the cores of the corresponding subgroups are  ${\rm PSL}_4(\Z_2) = {\rm
GL}_4(\Z_2) \cong \A_8$ for indices 8 and 15 and, for index 16, a semidirect product
$(\Z_2)^4  \rtimes  {\rm GL}_4(\Z_2)$ of order 322560; at present
we don't know if the abelianization of the core in the case of index 16 is
finite or infinite. On the other hand, the minimal indices of subgroups of
${\rm SOut} \, F_4$ are 8 and 15,  and index 16 is now missing.
We note that each ${\rm SAut} \, F_r$ admits a surjection onto a finite group
with a normal subgroup $(\Z_n)^r$ and factor group ${\rm GL}_r(\Bbb Z_n)$
(by dividing out first the kernel of the natural projection of
${\rm Inn} \, F_r  \cong  F_r$ onto $(\Z_n)^r$, then projecting onto
${\rm GL}_r(\Bbb Z_n)$).

\medskip

Employing in addition the quotient group procedure of GAP, we verified also
that the three smallest simple factor groups of ${\rm SAut} \, F_3$ and ${\rm
SOut} \, F_3$ are, as expected, the groups ${\rm PSL}_3(\Z_2)$ of order 168,
${\rm PSL}_3(\Z_3)$ of order 5616, and ${\rm PSL}_3(\Z_5)$ of order 372000 (we
note e.g. that, by the Lemma, the first Janko group does not occur since
it has no subgroup isomorphic to $SW_3 \cong S_4$).

\bigskip

We note that these computations can be slightly extended but that already the
suspected minimal index  31 of a proper subgroup of ${\rm SOut} \, F_5$ as well
as its suspected minimal nonabelian quotient ${\rm PSL}_5(\Z_2)$ (of order
9.999.360, with a subgroup of index 31) appear quite large for such
computations so we didn't pursue this further.

\end{document}